\documentclass[11pt]{article}
\usepackage[margin=1in]{geometry}
\usepackage{amsmath,amssymb,amsthm}
\usepackage{graphicx}
\usepackage{authblk}
\usepackage{float}
\usepackage{booktabs}
\usepackage{xcolor}
 \usepackage[numbers,sort&compress]{natbib} 
\usepackage{hyperref}
\usepackage{cleveref}

\title{\textbf{Economic resources, sporting efficiency and uncertainty: a global test  on the 2026 FIFA world cup}}

\author[1]{Thierno Thioune}
\affil[1]{Laboratory for Analysis, Research and Study of Development\\
University of Cheikh Anta Diop, Dakar, Senegal\\
\texttt{thierno.thioune@ucad.edu.sn}}

\author[2]{Babacar Mbaye Ndiaye}
\affil[2]{Laboratory of Mathematics of Decision and Numerical Analysis\\
University of Cheikh Anta Diop, Dakar, Senegal\\
\texttt{babacarm.ndiaye@ucad.edu.sn}}

\date{\today}

\begin{document}
\maketitle

\begin{abstract}
A related paper \citep{aathioune2026afcon} developed a stochastic optimization framework in which football performance combines GDP/population-based potential with a latent mean-reverting efficiency process.
The framework was applied to a cross-section of the 24 national teams participating in AFCON 2025.
This paper applies the same methodology to the full 48-team 2026 FIFA World Cup, extending the analysis beyond Africa without modifying the underlying model.
 Using GDP, population, FIFA ranking data, and complete tournament outcomes, it evaluates the framework against a global sample spanning all six confederations.
Two findings emerge. First, GDP and population explain only 9\% of single-tournament performance variance, compared with 33\% for FIFA points. This gap is consistent with the model’s distinction between systematic performance and match-specific idiosyncratic shocks.
Second, using the less noisy FIFA-points measure, the African pattern persists globally: Morocco, Senegal, Argentina, Croatia, Portugal, and Uruguay emerge as leading GDP-relative over-performers.
 Conversely, several major economies, including the United States, Germany, Japan, and South Korea, underperform their GDP-implied potential.
As in the related paper, we interpret the residual as a symmetric, non-structurally identified \emph{conditional performance residual}, rather than an SFA-consistent efficiency estimate. We also quantify the extent to which its cross-country variance mechanically reflects the underlying performance measure.
The GDP-rank/simulated-performance-rank correlation is robust to the choice of the model's un-identified dynamic parameters $(\kappa,\sigma)$ across a wide grid.
\end{abstract}


\section{Introduction}
The related paper on African football \citep{aathioune2026afcon} models national-team performance $Y_{i,t}$ as the product of an economic-potential production function $f(X_{i,t})$ and a latent, mean-reverting log-efficiency process $\theta_{i,t}=e^{z_{i,t}}$, plus an idiosyncratic shock $\varepsilon_{i,t}$:
\begin{equation}
\label{eq:model}
Y_{i,t} = \theta_{i,t} f(X_{i,t}) + \varepsilon_{i,t}, \qquad dz_{i,t} = \big[\kappa_i(\bar z_i - z_{i,t}) + \gamma_i u_{i,t}\big]dt + \sigma_i dW_{i,t}, \qquad \theta_{i,t}=e^{z_{i,t}}.
\end{equation}
The related paper derives \Cref{eq:model} as the state equation of a discounted stochastic control problem, gives the associated Hamilton--Jacobi--Bellman equation and optimal investment rule, and proves a verification theorem under explicit growth conditions \citep[Sections 3.2--3.5]{aathioune2026afcon}. It also positions the resulting efficiency measure explicitly relative to the stochastic frontier and data envelopment analysis (SFA/DEA) traditions \citep[Section 2.4]{aathioune2026afcon}. This paper does not repeat that derivation: it reuses only the state equation \eqref{eq:model} and the related paper's cross-sectional estimation strategy, applied here for the first time outside Africa. That paper estimated the model on a cross-section of the 24 Africa Cup of Nations 2025 teams. The 2026 FIFA World Cup, the first 48-team edition and jointly hosted by the United States, Canada, and Mexico, concluded on 19 July 2026 with Spain’s 1–0 extra-time victory over Argentina.
 It provides a natural large-scale out-of-sample test of whether the GDP/population–efficiency relationship observed for AFCON persists globally across all six confederations. \\
This paper applies the framework unchanged: the same production function $f(X_i) = A\cdot GDP_i^{\alpha_1} Pop_i^{\alpha_2}$, cross-sectional estimation with heteroskedasticity-robust inference, stochastic simulation of \eqref{eq:model}, and sensitivity analysis over $(\kappa,\sigma)$. 
As in \citep{aathioune2026afcon}, the collinear $GDPpc_i = GDP_i/Pop_i$ regressor is excluded. \\
As in the paper \citep{aathioune2026afcon}, we relabel the estimated residual a \emph{conditional performance residual} throughout, and we quantify explicitly how much of it is, by construction, a restatement of the performance measure used to estimate it (\Cref{sec:circularity}). \Cref{sec:sfa-dea} briefly restates, in this global setting, what such a residual is and is not relative to SFA/DEA. \Cref{sec:data} describes the data: all 48 qualified teams, their real GDP, population and FIFA ranking points, and the actual, complete 2026 World Cup result. \Cref{sec:results} reports two regressions with robust inference -- one on the actual tournament outcome, one on FIFA points -- compares them, and quantifies circularity for both. \Cref{sec:sensitivity} runs the simulation and sensitivity analysis. \Cref{sec:discussion} discusses what generalizes from Africa to the world, and what does not.

\section{Relation to stochastic frontier analysis and DEA}\label{sec:sfa-dea}
As discussed at length in the related paper \citep[Section 2.4]{aathioune2026afcon}, the residual estimated below is not a structurally identified efficiency term in the sense of stochastic frontier analysis \citep[SFA;][]{aigner1977} or data envelopment analysis \citep[DEA;][]{charnes1978}. Both SFA and DEA are static methods that decompose or bound a one-sided inefficiency term relative to a best-practice frontier, estimated from a snapshot of inputs and outputs. The quantity we estimate here by OLS, $z_i = \ln(Y_i) - \widehat{\ln f}(GDP_i,Pop_i)$, is a symmetric residual: it can be positive or negative, mixes ordinary noise with any genuine efficiency component, and no distributional assumption separates the two. We therefore refer to it throughout as a \emph{conditional performance residual} -- performance conditional on GDP and population -- rather than as a validated efficiency parameter. The \Cref{sec:circularity} quantifies exactly how much of its cross-country variance is mechanically determined by the fit of the underlying regression rather than by any independent efficiency signal. 
The contribution of both papers lies not in improving static efficiency estimation relative to SFA or DEA, but in providing a dynamic control framework once such an efficiency measure is observed. This distinction, discussed in \citep[Section 2.4]{aathioune2026afcon}, applies unchanged to the global sample considered here.
 
\section{Data}\label{sec:data}
\subsection{The 48 qualified teams}
The 2026 World Cup expanded the field from 32 to 48 teams for the first time. Slot allocation by confederation: UEFA 16, CAF 10 (9 direct plus DR Congo via the intercontinental play-off), AFC 9 (8 direct plus Iraq via the intercontinental play-off), CONCACAF 6 (three hosts -- Canada, Mexico, USA -- plus Panama, Haiti and Cura\c{c}ao), CONMEBOL 6, OFC 1 (New Zealand)\footnote{Compiled from FIFA/CAF/UEFA/CONCACAF qualification records and the March 2026 inter-confederation play-off tournament (winners: DR Congo and Iraq, from a six-team pool that also included Bolivia, Jamaica, New Caledonia and Suriname).}. 
Economic data: nominal GDP and GDP per capita, from the IMF World Economic Outlook database directly (April 2026 edition), 2026 projections;\footnote{International Monetary Fund, World Economic Outlook database, April 2026 edition, \url{https://www.imf.org/en/Publications/WEO}.} population is computed as $GDP/GDPpc$. 
FIFA ranking data: FIFA Men's World Ranking points, from the FIFA World Ranking directly, ranking dated 20 July 2026.\footnote{FIFA World Ranking, \url{https://www.fifa.com/fifa-world-ranking}, ranking date 20 July 2026.}

\noindent \textbf{A data limitation, disclosed.} England and Scotland compete as separate national teams but are not sovereign states and have no independent GDP or population statistics; both are assigned the United Kingdom's GDP and population as the best available proxy. \\
This is a genuine approximation, not a modeling choice. It mechanically produces a large negative conditional performance residual for whichever of the two has the lower FIFA ranking (Scotland), simply because both are compared to the same (larger, UK-wide) economic potential. 
We flag every instance where this affects the results below.

\subsection{The actual 2026 World Cup result}
Champion: Spain. Runner-up: Argentina. Third place (beating France 6--4 in a ten-goal match): England. Fourth place: France. Quarter-final losers: Morocco (to France), Belgium (to Spain), Norway (to England), Switzerland (to Argentina).\footnote{Cross-checked against contemporaneous match reports (NPR, CBS News, ESPN, beIN Sports; July 2026).} Following the coarse, independently verifiable scoring scale of \citep{aathioune2026afcon} (chosen for direct comparability): Champion$=6$, Runner-up$=5$, 3rd$=4$, 4th$=3.5$, quarter-final exit$=2.5$, and every other team (round of 32, round of 16 or group-stage exit)$=1$.

\section{Results}\label{sec:results}
\subsection{Two regressions, robust inference, and why they disagree}
We estimate, by OLS on the $N=48$ qualified teams,
\[
\ln(Y_i) = \ln A + \alpha_1 \ln GDP_i + \alpha_2 \ln Pop_i + z_i,
\]
once with $Y_i$ equal to the actual World Cup 2026 performance score (regression A), and once with $Y_i$ equal to FIFA ranking points (regression B, our primary measure). \Cref{tab:regressions} reports both, with heteroskedasticity-robust (HC3) standard errors and the joint test $H_0:\alpha_1=\alpha_2=0$, appropriate given the modest sample size and the very different variance structure of the two outcome variables.

\begin{figure}[h!]
\centering
\includegraphics[width=0.55\linewidth]{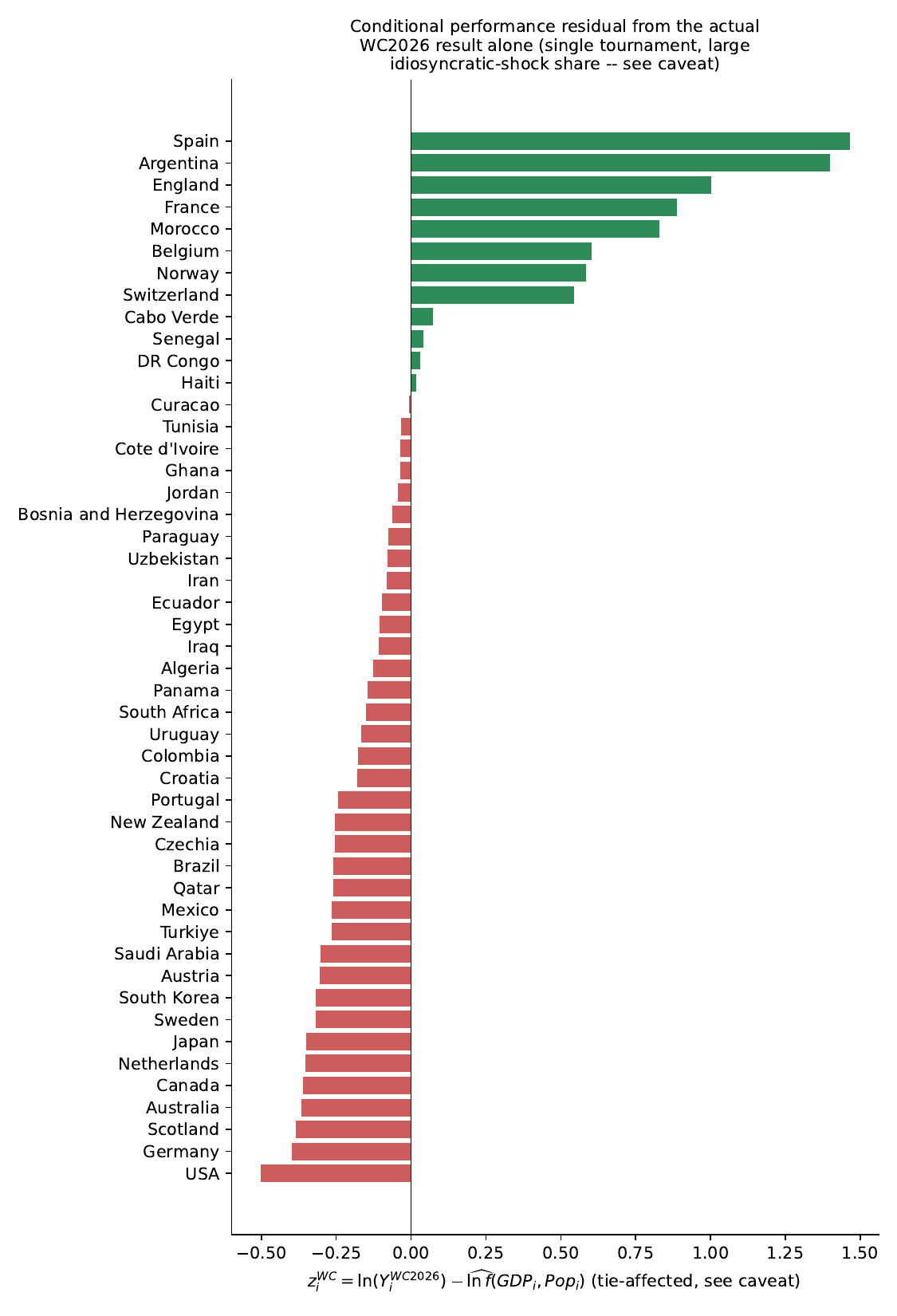}
\caption{Conditional performance residual from the actual WC2026 result alone (regression A): 40 of 48 teams are tied, so their ordering is driven almost mechanically by $-\widehat{\ln f}(GDP_i,Pop_i)$.}
\label{fig:wcresidual}
\end{figure}

\begin{table}[h!]
{\footnotesize
\centering
\caption{Two cross-sectional regressions on the same 48 teams. HC3 robust standard errors in parentheses.}
\label{tab:regressions}
\begin{tabular}{lrrrl}
\toprule
 & $\hat\alpha_1$ (GDP) & $\hat\alpha_2$ (Population) & $R^2$ (adj.) & Joint test $H_0:\alpha_1=\alpha_2=0$ \\
\midrule
(A) Actual WC2026 score & 0.109 (0.050)$^{**}$ & $-0.065$ (0.047) & 0.091 (0.050) & $F=2.44$, $p=0.099$ \\
(B) FIFA points \emph{(primary)} & 0.034 (0.011)$^{***}$ & 0.003 (0.012) & 0.328 (0.298) & $F=13.32$, $p<0.001$ \\
\bottomrule
\end{tabular}
\vspace{2pt}
\\ \footnotesize $^{**}p<0.05$, $^{***}p<0.01$ (robust $t$-tests, individual coefficients).
}
\end{table}
\noindent The gap between the two $R^2$s (0.09 vs.\ 0.33) is itself a substantive finding, not noise to be explained away, and robust inference sharpens rather than changes this reading. 
In regression (A), 40 of the 48 teams share the same pre-quarter-final score, making their residuals mechanically equal to $z_i = -\widehat{\ln f}(GDP_i,Pop_i)$.
Thus, these residuals largely reflect the coarse scoring structure rather than an independent performance signal (\Cref{fig:wcresidual}), explaining the extreme negative residual of the United States.
 
\noindent Although GDP is individually significant in regression (A) (($p=0.029$), the joint test is only marginally significant (($p=0.099$) and explanatory power remains low. With a coarse, heavily tied outcome, these results are not contradictory and should be interpreted jointly rather than in isolation.
More fundamentally, a single-elimination tournament provides only a thin realization of the idiosyncratic shock $\varepsilon_{i,t}$, whereas FIFA points aggregate performance over a much longer horizon. By averaging out these shocks, FIFA points better capture the systematic component $\theta_{i,t}f(X_{i,t})$ targeted by the model. 
This is exactly what the related paper's own state-space formulation anticipates: a noisier, higher-frequency $Y$ identifies $\theta$ less precisely than a smoothed one. 
We therefore use regression (B) -- FIFA points, directly comparable to the AFCON paper's own measure, and jointly significant at $p<0.001$ -- as the primary measure for the simulation below, and use (A) only as a qualitative check. \\ 
\textbf{A validation check.} Despite the low $R^2$ of regression (A), the FIFA-implied performance ranking exhibits predictive content for actual tournament outcomes. Of the eight 2026 quarter-finalists, six (Argentina, Spain, Morocco, France, England, and Belgium) rank among the FIFA-implied top 10 by $z_i$. FIFA World Cup 2026 Quarterfinal Fixtures confirms these six were among the eight quarter-finalists. 
Norway and Switzerland, the exceptions, are precisely two of the tournament's best-documented surprise runs.

\subsection{Quantifying circularity in both regressions}\label{sec:circularity}
As emphasized in the revised related paper, the correlation between an OLS residual and its dependent variable is not free evidence of an independent signal: by construction, $\mathrm{Corr}(z_i,\ln Y_i) = \sqrt{1-R^2}$ for \emph{any} regression, regardless of what $Y_i$ measures. We compute this explicitly for both regressions here. \\
For regression (A) (actual WC2026 score, $R^2=0.091$),
\[
\mathrm{Corr}\big(z_i^{WC},\, \ln Y_i^{WC2026}\big) = \sqrt{1-0.091} = 0.954,
\]
confirmed numerically at $0.9537$. This extremely high mechanical correlation is the same phenomenon already flagged qualitatively above (the tie artifact): with such a low $R^2$, the residual is, for practical purposes, almost the raw (log) outcome itself, re-centered. \\
For regression (B) (FIFA points, $R^2=0.328$),
\[
\mathrm{Corr}\big(z_i,\, \ln(\mathrm{FIFA}_i)\big) = \sqrt{1-0.328} = 0.820,
\]
confirmed numerically at $0.8197$, and the Spearman rank correlation between the country ranking implied by $z_i$ and the raw FIFA-points ranking is $0.801$ ($p<0.001$). This is close to, and slightly higher than, the corresponding figure in the AFCON study ($0.808$ for the Pearson correlation, $0.734$ for the Spearman rank correlation, on $N=24$). \\
As in the related paper, claims of GDP-relative over-performance should be interpreted in light of the modest explanatory power of GDP and population ($R^2=0.33$).
Thus, ``over-performance'' denotes a conditional deviation from the model-implied benchmark rather than a direct measure of efficiency.\\
The residual likely also captures omitted factors, including institutional quality, diaspora and expatriate players, football tradition, and demographic structure. Reliable data on these factors are unavailable across the full 48-country, six-confederation sample, limiting their explicit inclusion in the analysis.
We requalify $z_i$ as a \emph{conditional} performance residual accordingly, including in the abstract, rather than as a validated efficiency parameter.

\subsection{Conditional performance residuals (FIFA points, primary measure)}
\begin{figure}[h!]
\centering
\includegraphics[width=0.55\linewidth]{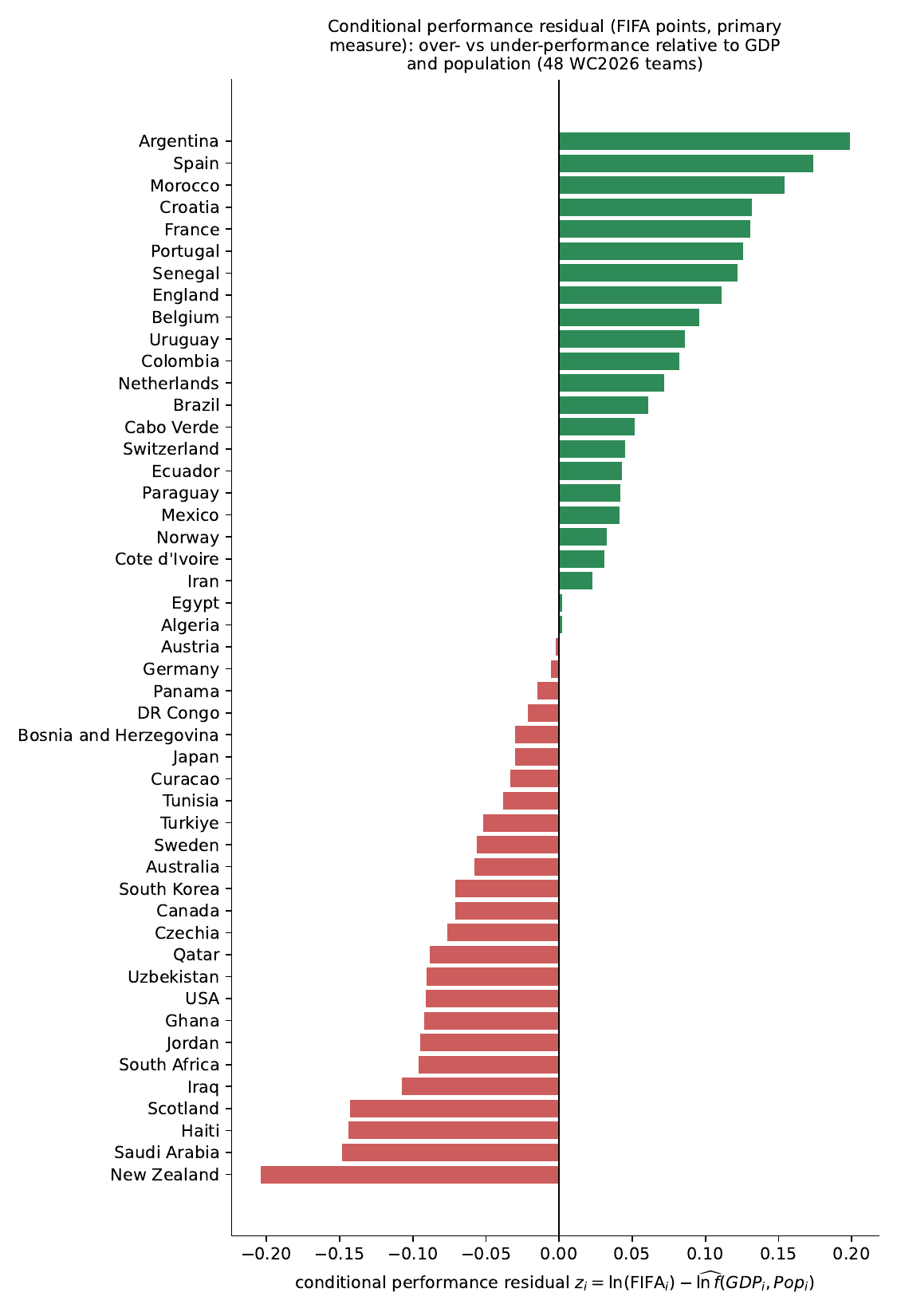}
\caption{Conditional performance residual $z_i$ (FIFA points regression, primary measure) for all 48 World Cup 2026 teams.}
\label{fig:residual}
\end{figure}

\Cref{fig:residual} plots $z_i$ for all 48 teams. 
The largest positive residuals are Argentina, Spain, Morocco, Croatia, France, Portugal, Senegal, England, Belgium and Uruguay. 
The largest negative residuals are New Zealand, Saudi Arabia, Haiti, Scotland (data-limitation caveat above), Iraq, South Africa, Jordan, Ghana, USA and Uzbekistan. 
Morocco and Senegal, the two leading African over-performers in \citep{aathioune2026afcon}, again rank among the largest over-performers globally. This provides out-of-sample support for the related paper’s central finding, albeit subject to the circularity caveat discussed above.\\
 Several of the world's largest economies under-perform their GDP/population potential: the USA (\$32.4 trillion GDP, largest in the field, $z=-0.09$), Germany ($z\approx -0.01$, borderline), Japan ($z=-0.03$), South Korea ($z=-0.07$) and Canada ($z=-0.07$).

\section{Simulation and sensitivity}\label{sec:sensitivity}
\subsection{Simulation and GDP-rank vs.\ performance-rank}
\begin{figure}[h!]
\centering
\includegraphics[width=0.6\linewidth]{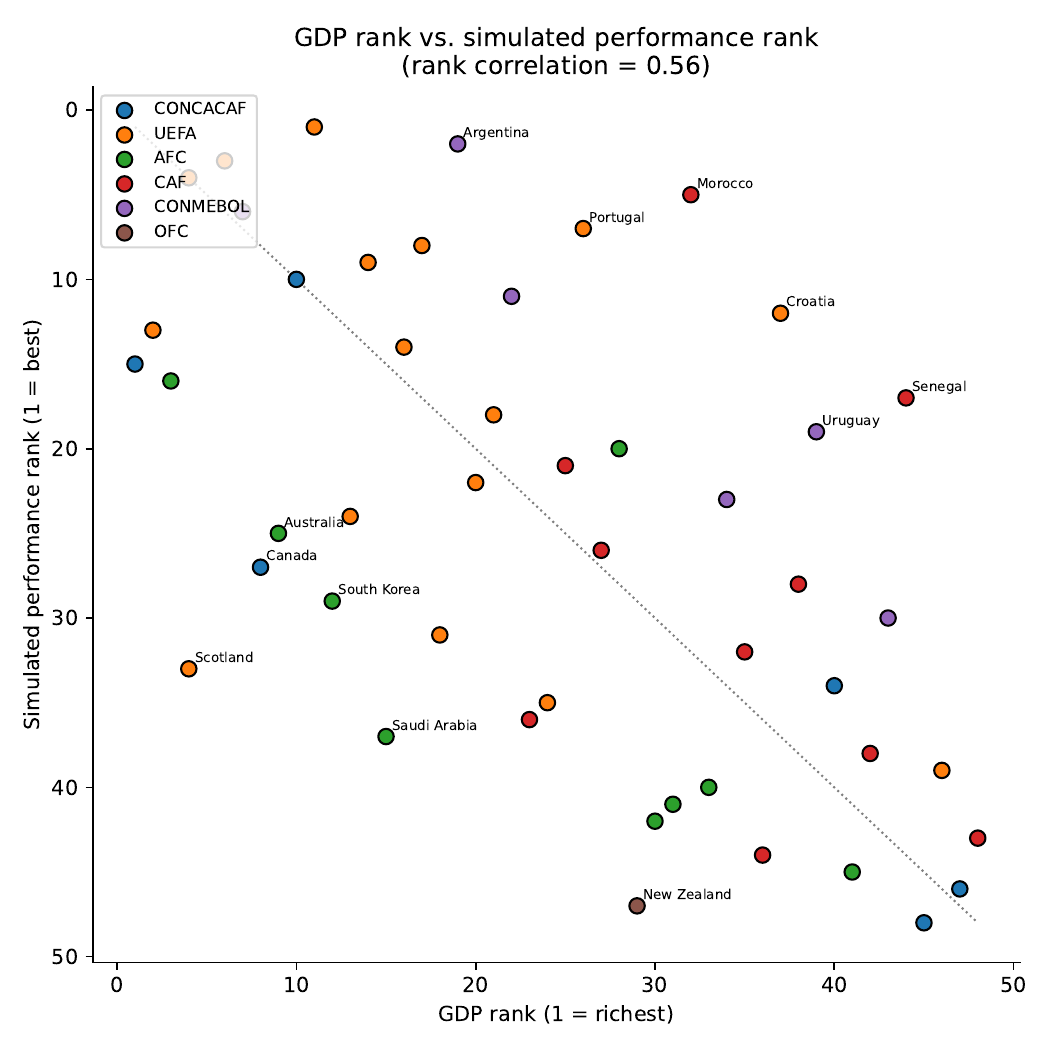}
\caption{GDP rank versus simulated performance rank, 48 World Cup 2026 teams, by confederation.}
\label{fig:rankscatter}
\end{figure}

As in \citep{aathioune2026afcon}, each team is simulated forward from $z_{i,0}=\bar z_i = z_i$ under \eqref{eq:model} with shared illustrative parameters ($\kappa=0.7$, $\sigma=0.30$, $\gamma=0.6$, $c=1$, constant investment $u=0.2$, horizon $T=5$ years, $f_i$ the GDP/population-based fitted potential). \\
As in the related paper, setting $\bar z_i=z_i$ is an identification assumption rather than a normalization. It enables simulation from a single cross-section without estimating $\bar z_i$, but mechanically preserves each country’s initial relative standing rather than identifying it endogenously. \\ 
\Cref{fig:rankscatter} shows a positive rank correlation ($0.562$) between GDP and simulated performance, while substantial deviations highlight the role of the conditional performance residual.
 This pattern is consistent with the AFCON analysis, where the corresponding correlation was $0.62$.
 
\subsection{Sensitivity to $(\kappa,\sigma)$}
\begin{figure}[h!]
\centering
\includegraphics[width=0.6\linewidth]{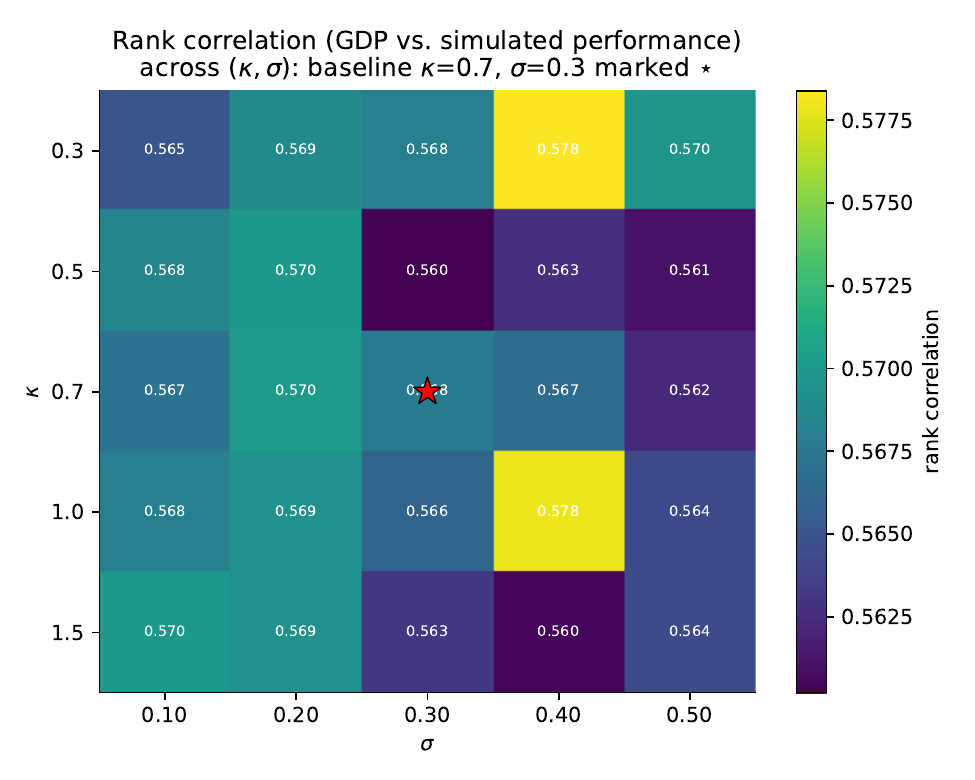}
\caption{Rank correlation between GDP and simulated performance across a $5\times 5$ grid of $(\kappa,\sigma)$, 48 World Cup 2026 teams; baseline marked with a star.}
\label{fig:heatmap}
\end{figure}

As in the AFCON study, $\kappa$ and $\sigma$ are not calibrated, since a single cross-section cannot identify dynamic parameters. The same limitation arose in \citep[Section 6.1.1]{aathioune2026afcon}, where calibration from a four-edition AFCON panel was abandoned rather than reported as a potentially spurious point estimate. 
We check robustness over the same grid, $\kappa\in\{0.3,0.5,0.7,1.0,1.5\}$, $\sigma\in\{0.10,0.20,0.30,0.40,0.50\}$ (Monte Carlo, $N=8{,}000$ paths per cell). 
The rank correlation stays in a narrow band, $[0.560,0.578]$, across all 25 combinations (\Cref{fig:heatmap}) -- comparable to, if slightly wider than, the AFCON result ($[0.605,0.643]$ under the same recomputation). \\
Across the parameter grid, the five highest-ranked teams take only two closely related configurations: \{Argentina, Brazil, England, France, Spain\} or \{Argentina, England, France, Morocco, Spain\}.
 The five lowest-ranked teams are similarly stable, with the two configurations differing only by Cabo Verde and Ghana.\\
At the baseline parameters, Morocco, Senegal, Croatia, Uruguay, and Portugal are the five largest  \emph{over}-performers relative to GDP rank.
 Conversely, Scotland, Saudi Arabia, Canada, New Zealand, and South Korea exhibit the largest  \emph{under}-performance. \\
The central empirical pattern persists from the 24-team AFCON sample to the 48-team global field: GDP matters, while the conditional performance residual has a stronger marginal effect.
This result is robust to the choice of the unidentified structural parameters $\kappa$ and $\sigma$.
 
\section{Discussion}\label{sec:discussion}
\subsection{What generalizes from Africa to the world}
Three results from \citep{aathioune2026afcon} replicate at global scale. First, GDP and population alone leave the majority of cross-country variation in footballing strength unexplained ($R^2=0.33$ here vs.\ $0.35$ in the 24-country AFCON regression) -- essentially the same order of magnitude, now confirmed on a sample twice the size and drawn from every confederation, not just one. \\
In both cases, roughly two-thirds of the residual’s cross-country variance mechanically reflects the underlying performance measure (\Cref{sec:circularity}). Thus, this consistency indicates stability in the model’s fit rather than independent evidence of an efficiency channel.\\
Second, Morocco and Senegal emerge as leading over-performers in both studies. This provides out-of-sample corroboration, as the World Cup analysis uses a distinct comparison set of 46 non-African teams and an independent tournament outcome. \\
Third, the GDP-rank/performance-rank correlation is positive but firmly below 1 in both settings ($0.562$ here, $0.62$ in Africa), and in both cases this finding is robust to the choice of the model's un-identified dynamic parameters over a wide grid.

\subsection{What is new here}
The World Cup exercise also surfaces a finding the single-confederation AFCON study could not: the sharp contrast between a single-tournament outcome and an aggregated ranking as the performance measure $Y_i$. 
With $R^2=0.09$, World Cup outcomes are only weakly explained by GDP and population, while the tied scoring of 40 mid-table exits induces a strong mechanical residual–outcome correlation ($0.954$; \Cref{sec:circularity}).
By contrast, FIFA points, aggregated over several years and many matches, are substantially more explained by economic fundamentals ($R^2=0.33$), reinforcing the GDP/conditional-performance distinction. \\
This is a clean illustration, with real data, of the model's own distinction between the systematic term $\theta_{i,t}f(X_{i,t})$ and the idiosyncratic shock $\varepsilon_{i,t}$: single knockout tournaments are dominated by the latter, long-run rankings by the former. \\ 
This finding cautions against interpreting any single tournament, African or global, as a clean measure of genuine efficiency.  It reinforces the need for a multi-edition panel to calibrate $\kappa$ and ,$\sigma$ and to identify the conditional performance residual using richer controls and a properly specified SFA decomposition.
 
\subsection{Remarks}
As in \citep{aathioune2026afcon}: this is a single cross-section, not a panel, so $z_i$ is a snapshot (``footballing strength relative to GDP and population, as of the concluded World Cup 2026''), not a fully identified structural efficiency parameter in the SFA sense (\Cref{sec:sfa-dea}), and $\kappa,\gamma,c,\sigma$ are shared across all 48 teams by necessity, with the $\bar z_i=z_i$ identification assumption flagged explicitly in \Cref{sec:sensitivity}. The England/Scotland UK-proxy issue (\Cref{sec:data}) mechanically penalizes Scotland's residual and should not be read as a substantive claim about Scottish football governance. Finally, the actual-result regression's low $R^2$ should not be read as ``GDP does not matter for World Cup performance''. 
This reflects the small-sample nature of a single knockout tournament. The FIFA-points regression, based on the same 48 teams and covariates and jointly significant at $p<0.001$, nevertheless reproduces the GDP/conditional-performance split observed in Africa.

\section{Conclusion}
Applying the framework of \citep{aathioune2026afcon} unchanged to the 48-team 2026 FIFA World Cup confirms its central finding at the global level: economic resources determine potential, while a sizeable GDP-independent component drives performance above or below it.
Following the related paper’s revised methodology, we interpret this component as a conditional residual, noting its mechanically strong correlation with the underlying performance measure ($\mathrm{Corr}=0.82$) and its absorption of unmodeled factors.
Morocco and Senegal are over-performers both within Africa and globally; several of the world's richest football economies under-perform their potential, with robust statistical support for the joint explanatory role of GDP and population ($p<0.001$) despite individually modest elasticities. 

\bibliographystyle{unsrtnat}   
\bibliography{references}       

@article{aathioune2026afcon,
      author={Thierno Thioune and Babacar Mbaye Ndiaye},
      title={Economic resources, sporting efficiency and uncertainty: a stochastic optimization model of football performance in Africa}, 
      year={2026},
      eprint={2610.06967},
      archivePrefix={arXiv},
      primaryClass={math.OC},
      url={https://arxiv.org/abs/2610.06967}, 
}

@article{aigner1977,
  author  = {Aigner, Dennis and Lovell, C. A. Knox and Schmidt, Peter},
  title   = {Formulation and Estimation of Stochastic Frontier Production Function Models},
  journal = {Journal of Econometrics},
  year    = {1977},
  volume  = {6},
  number  = {1},
  pages   = {21--37}
}

@article{charnes1978,
  author  = {Charnes, Abraham and Cooper, William W. and Rhodes, Edwardo},
  title   = {Measuring the Efficiency of Decision Making Units},
  journal = {European Journal of Operational Research},
  year    = {1978},
  volume  = {2},
  number  = {6},
  pages   = {429--444}
}

\end{document}